\documentclass[11pt]{amsart}
\usepackage{amsmath,amssymb,mathrsfs,amscd}
\usepackage{graphicx}

\title{Twistor lines on Nagata threefold}

\author{Nobuhiro Honda$^{\dag}$}
\thanks
{$^{\dag}$This work was partially supported by
Research Fellowships of the 
Japan Society for the Promotion
of Science for Young Scientists.
}
\date{}

\newcommand{\ol}{\overline}
\newcommand{\ra}{\rightarrow}
\newcommand{\lra}{\longrightarrow}

\newcommand{\set}{\,|\,}
\newcommand{\proofend}{\hfill$\square$}

\newcommand{\bsl}{\backslash}

\newtheorem{prop}{Proposition}[section]
\newtheorem{lemma}[prop]{Lemma}

\newtheorem{thm}[prop]{Theorem}
\newtheorem{rmk}[prop]{Remark}

\begin{document}
\maketitle
\begin{abstract}
We give an explicit  description of rational curves in the product of three copies of complex projective lines, which are transformed into twistor lines in M. Nagata's example of non-projective complete algebraic variety,
viewed as the twistor space of Eguchi-Hanson metric. 
In particular, we show that there exist two families of such curves and both of them are parameterized by mutually diffeomorphic, connected real 4-dimensional manifolds.
We also give a relationship between these two families through  a birational transformation naturally associated to the Nagata's example.
\end{abstract}

\bigskip\noindent
\section{ Introduction.}
In 1958, M. Nagata \cite{Na58} constructed a remarkable example of a compact complex threefold, which was the first example of a non-projective complete algebraic variety.
A new aspect of this threefold was given by A. Fujiki \cite{F}, who showed  that  Nagata's example is a compactification of the twistor space of the famous Eguchi-Hanson metric on the cotangent bundle of $\mathbf P^1$ \cite{EH78}.
Fujiki proved this result by investigating the Calabi family naturally associated to a hyper-K\"ahler metric.
The most important geometric object in a twistor space is the twistor lines.
Since Nagata's example is obtained from a product $\mathbf P^1\times\mathbf P^1\times\mathbf P^1$ by applying a simple birational transformation, it seems natural to ask which curves in $\mathbf P^1\times\mathbf P^1\times\mathbf P^1$ are transformed into  twistor lines in Nagata's example.
In this note we  give an answer to this question by describing all of them explicitly.

After recalling  Nagata's construction,
we first determine which curves in $\mathbf P^1\times\mathbf P^1\times\mathbf P^1$ can be transformed into  real lines in the Nagata threefold  (Propositions \ref{prop-type} and \ref{prop-explicit}).
Here, a smooth rational curve in a threefold is called a {real line} if the normal bundle of the curve is isomorphic to $\mathscr O(1)\oplus \mathscr O(1)$ and  it is invariant under an anti-holomorphic involution of the threefold.
These are candidates of twistor lines.
The parameter space of these candidates consists of two connected 4-manifolds, which will be denoted by $M^+$ and $M^-$.
There is a natural isomorphism between $M^+$ and $M^-$ which comes from an involution of $\mathbf P^1\times\mathbf P^1\times\mathbf P^1$ of the same degree.
Although $M^+$ and $M^-$ parameterize real curves in $\mathbf P^1\times\mathbf P^1\times\mathbf P^1$, 
both contain 2-dimensional submanifolds parameterizing reducible curves.
We find that, as a consequence of Nagata's birational transformation,
this situation is resolved for precisely one of $M^+$ and $M^-$ so that  all  its curves become irreducible.
Then by using our explicit description of the candidates, we show that the resolved  family actually becomes the set of twistor lines (Theorem \ref{thm-main}).

All irreducible curves of the other family are also transformed into real lines in the Nagata threefold.
But the family cannot be the set of twistor lines since the reducible curves (parameterized by the 2-dimensional submanifold mentioned above) remain reducible in Nagata threefold (even after applying the birational transformation).
We observe that all these reducible curves contain a common curve (which will be written by $B_0^++B_{\infty}^+$)  which is homologous to zero.
We will also find that  any members of this family can be deformed into a twistor line (keeping the reality and smoothness) if one allows the members to pass the situation that they have $\mathscr O(2)\oplus \mathscr O$ as their normal bundles in the Nagata threefold.
In particular,  all irreducible curves in the other family are homologous to actual twistor lines.


Finally, we mention that in the paper \cite{Hi79} N.\,J. Hitchin gave an explicit description of the twistor spaces of multi-Eguchi-Hanson metrics (i.\,e. Gibbons-Hawking metrics) and their twistor lines.
For the twistor spaces of multi-Eguchi-Hanson metrics, explicit birational transformation into a rational variety of a simple form seems unknown, except for the original Eguchi-Hanson metric, which is  adapted in this papar.
Also, in the paper \cite{B86}, D. Burns explicitly constructed the twistor spaces of the contangent bundles of Hermitian symmetric spaces and their twistor lines.
The author does not know how to identify  these works with the explicit description in this paper.

\bigskip

\noindent{\bf Acknowledgment}
The author would like to thank the referee for careful reading and invaluable comments. 
Also he would like to thank Shihoko Ishii for answering his question on a Hilbert scheme.

\section{Explicit description of real lines}
First of all, we recall the construction of Nagata threefold \cite{Na58}.
We write $\mathbf P=\mathbf P^1$ for the complex projective line throughout this paper.
Let $X=\mathbf{P}\times\mathbf{P}\times\mathbf{P}$ be the product of three complex projective lines, and write $C=\mathbf{P}$ for the last factor  which  plays a special role.
Let  $f:X\ra C$ be  the projection, and
 $0$ and $\infty$  two distinct points of $C$.
We fix an isomorphism between the first and the second factors of $X$ and let $\Delta\subset\mathbf{P}\times\mathbf{P}$ be the diagonal.
Put $Q=\Delta\times C$, a divisor on $X$.
For $t\in C$, we write $\Delta_t$ for  the diagonal of $X_t:=f^{-1}(t)$.
Let $\mu:Y\ra X$ be the blowing-up  along $\Delta_0\cup \Delta_{\infty}$, and $X_0'$ and $X'_{\infty}$ the strict transforms of $X_0$ and $X_{\infty}$ respectively.
Because $\Delta_0\subset X_0$ and $ \Delta_{\infty}\subset  X_{\infty}$, $\mu$ gives isomorphisms $X'_0\simeq X_0$ and $X'_{\infty}\simeq X_{\infty}$. 
Moreover, $N_{X'_0/Y}$ and $N_{X'_{\infty}/Y}$, the normal bundles of $X'_0$ and $X'_{\infty}$ in $Y$ respectively, are isomorphic to $\mathscr O(-1,-1)$, where $\mathscr O(m,n)$ denotes the holomorphic line bundle on $\mathbf{P}\times\mathbf{P}$ of bidegree $(m,n)$.
Therefore  $X'_0\, (\simeq\mathbf P\times\mathbf P)$  can be blown-down to $\mathbf{P}$ along each two projections, and 
 the same thing holds  for $X'_{\infty}$. 
Let $\nu_0^+:X'_0\ra\mathbf{P}$ be  the projection to the {\em  first} factor and $\nu_{\infty}^+:X'_{\infty}\ra \mathbf{P}$ the projection to the {\em  second} factor.
Let $\nu^+:Y\ra Z^+$ be the blowing down of $Y$ inducing $\nu_0^+$ and $\nu_{\infty}^+$ on $X'_0$ and $X'_{\infty}$ respectively.
If we interchange  the role of the two factors, 
we obtain another threefold $Z^-$.
In this way we obtain two smooth threefolds $Z^+$ and $Z^-$. 
Both $Z^+$ and $Z^-$ are called {\em{Nagata threefolds}}. 
The isomorphism of $X$ interchanging the first and the second factors naturally induces a biholomorphic map $i:Z^+\ra Z^-$ and hence these are naturally biholomorphic.
We denote by $f_+:Z^+\ra C$   for the projection naturally induced from $f:X\ra C$.
This is a  holomorphic submersion.
The fibers $f_+^{-1}(0)$ and $f_+^{-1}(\infty)$ are  isomorphic to $\Sigma_2=\mathbf{P}(\mathscr O\oplus \mathscr O(2))$ and all the other fibers are isomorphic to $\mathbf P\times\mathbf P$.
If we denote $B_0^+=\nu^+(X'_0)$ and $B^+_{\infty}=\nu^+(X'_{\infty})$, $B_0^+$ and $B^+_{\infty}$ are the minimal section of $\Sigma_2\ra\mathbf{P}$.
It was shown in \cite{Na58} that the effective curve $B_0^++B^+_{\infty}$ is homologous to zero and hence $Z^+$ does not admit a K\"ahler metric.
Finally, we define a divisors $Q^+$ in $Z^+$  as the bimeromorphic image of  $Q=\Delta\times C$ into $Z^+$.
Similarly $f_-:Z^-\ra C$, $B_0^-, B_{\infty}^-$ and $Q^-$ are defined for $Z^-$.
We have $i(Q^+)=Q^-$.

Next we introduce real structures.
Let $\sigma_1$ be an anti-podal map of $\mathbf{P}$, which is an anti-holomorphic involution without fixed points.
$\sigma_1$ can be explicitly given by $\sigma_1(z_0:z_1)=(-\ol{z}_1:\ol{z}_0)$ using homogeneous coordinates.
Then we define an anti-holomorphic involution $\sigma$ on $X$ by
\begin{equation}\label{eqn-rs}\sigma(x,y,z)=(\sigma_1(y),\sigma_1(x),\sigma_1(z)).
\end{equation}
Here we are using the identification of the first and the second factors.
Clearly $\sigma$ has no fixed points.
Further $\sigma$ survives under the above construction and induces  real structures on the Nagata threefolds. 
We  denote these real structures by $\sigma_+$ on $Z^+$ and $\sigma_-$ on $Z^-$ respectively. 
Then $Q^{\pm}\subset Z^{\pm}$ are invariant under $\sigma_{\pm}$; namely they are real.
Moreover, it is immediate to see $\sigma_-i=i\sigma_+$ and hence $(Z^+,\sigma_+)$ and $(Z^-,\sigma_-)$ are isomorphic  as complex manifolds with real structures.
So in the sequel we mainly consider $(Z^+,\sigma_+)$;
corresponding results for $(Z^-,\sigma_-)$ are immediately obtained by using the isomorphism $i$.

We are going to show that $(Z^+\backslash Q^+,\sigma_+)$ is the twistor space of a hyperK\"ahler metric on the total space of the cotangent bundle of  $\mathbf P$, mainly by showing that the complex threefold $Z^+\backslash Q^+$ is actually foliated by real lines.
As a first step, we show the following proposition  describing the images of real lines in $Z^+$ (if any) under the above birational correspondence between  $Z^+$ and $ X$.
To indicate curves in $X$, we use the following notation: 
by taking three factors of $X$ as generators,
we have a natural isomorphism $H_2(X,\mathbf Z)\simeq \mathbf Z^3$, and the homology class of an effective curve on $X$ is determined by  $(k_1,k_2,k_3)\in\mathbf Z^3$;
in other words,  $k_i$ is the intersection number of the curve with  fibers of the $i$-th projection from $X$ to $\mathbf P$.
We call $(k_1,k_2,k_3)$ {\em{the degree}}\, of a curve for simplicity.
On the other hand, a {\em real line} in a threefold (with real structure)  means a smooth rational curve which is invariant under the real structure and which has $\mathscr O(1)\oplus\mathscr O(1)$ as the normal bundle.

The following proposition gives a necessary condition for a rational curve in $X=\mathbf P\times\mathbf P\times\mathbf P$ to be transformed into a real line in the Nagata threefold.

\begin{prop}\label{prop-type}
Let $(Z^+,\sigma_+)$ be the Nagata threefold equipped with the real structure given above, and suppose that $L$ is a real line in $Z^+$ which is disjoint from the divisor $Q^+$.
Then $L$ must be a section of $f_+\!:\!Z^+\ra C$.
Further by the birational correspondence between  $Z^+$ and $ X$, $L$ is transformed into a real curve in $X$ satisfying the following.
(i) If $L$ does not go through $B_0^+$, then the degree of the corresponding  curve in $X$ is $(1,1,1)$ and it goes through $\Delta_0$ (and hence $\Delta_{\infty}$ also).
(ii) If $L$ goes through $B_0^+$ (and hence $B_{\infty}^+$ also), then the degree of the corresponding curve in $X$ is $(0,0,1)$.
\end{prop}

\noindent Proof.
Assume $L$ is a real line disjoint from $Q^+$.
First we show that $L$ is a section of $f_+$.
By standard calculations, we can show that the anti-canonical line bundle of $Z^+$ satisfies
\begin{equation}\label{eqn-can}
-K_{Z^+}\simeq  \mathscr O(2Q^+) \otimes f_+^*\mathscr O(4).
\end{equation}
On the other hand, by adjunction formula, we have $-K_{Z^+}\cdot L=4$.
Hence since we have assumed that $L$ is disjoint from $Q^+$, we have $f_+^*\mathscr O (1)\cdot L=1$ and thus $L$ is a section of $f_+$.

Next to prove (i) suppose that $L$ is a real line   in $Z^+$ that is disjoint from $B_0$ (and $B_{\infty})$, and let $\nu^+:Y\ra Z^+$ and $\mu:Y\ra X$ be as in the construction of $Z^+$.
$L$ is a section of $f_+$.
Obviously $\nu^+$ does not change a neighborhood of $L$ in $Z$ and therefore $L$ still intersects $X'_0$ and $X'_{\infty}$ once respectively in $Y$.
Hence  under the blowing-down $\mu:Y\ra X$ the image of the real line goes through $\Delta_0=\mu(X'_0)$ and $\Delta_{\infty}=\mu(X'_{\infty})$.
Next under the same assumption we show that the degree of the image of $L$ in $X$ is $(1,1,1)$.
Because the image  remains to be a section of $f:X\ra C$, the degree must be of the form $(k_1,k_2,1)$ for non-negative integers $k_1$ and $k_2$.
Moreover, since the birational transformations keep the reality of $L$, and since the real structure on $X$ interchanges the first and second factors as in  (\ref{eqn-rs}), we have $k_1=k_2$.
Moreover, the degree of the normal bundle must increase by two under the blowing-down $\mu:Y\ra X$. 
Therefore it becomes $2+2=4$. 
It is readily seen that this  happens iff $k_1=(k_2=)1$. 
Thus we obtain (i) of the proposition.

Next suppose that $L$ is a real line  intersecting $B_0^+$ (and $B_{\infty}^+)$.
Then since $L$ is  a section of $f_+$, 
$L$ intersects $f_+^{-1}(0)$ and $f_+^{-1}(\infty)$ transversally at a unique point respectively.
Hence the degree of the normal bundle in $Y$ becomes $2-2=0$.
Moreover, this time, the blowing down $\mu:Y\ra X$ makes no effect in a neighborhood of the curve.
Hence the degree of the normal bundle of the image of $L$ in $X$ must be zero.
This happens only when  the degree is $(0,0,1)$ and we obtain (ii) of the proposition.
\proofend

\vspace{3mm}
The real curves in (i) and (ii) of Proposition \ref{prop-type} can be easily written down explicitly.
Let $(x,y,t)$ be an affine coordinate on $X$ so that the real structure $\sigma$ is given by  $\sigma(x,y,t)=(-1/\ol{y},-1/\ol{x},-1/\ol{t})$.
Then we have the following.

\begin{prop}\label{prop-explicit}
(i) A real irreducible curve  of degree $(1,1,1)$ in $X$ going through the point $(d,d,0)\in\Delta_0$ with $d\in\mathbf C$ is of the form
\begin{equation}\label{eqn-111}
x=\frac{d-at}{1+a\ol{d}t},\,\,y=\frac{\ol{a}d-t}{\ol{a}+\ol{d}t},
\end{equation}
for some  $a\in \mathbf C^*$, whereas if $d=\infty$ it is of the form
\begin{equation}\label{eqn-111'}
x=\frac{1}{at},\,\,\,y=\frac{\ol{a}}{t},
\end{equation}
for some $a\in\mathbf C^*$.
(ii) A real curve in $X$ of degree $(0,0,1)$ is given by 
\begin{equation}\label{eqn-001}
x=d,\,\,y=-\frac{1}{\ol{d}},
\end{equation}
where $d\in\mathbf C\cup\{\infty\}$.
(iii) If $|a|=1$, the curves \eqref{eqn-111} and \eqref{eqn-111'} are contained in $Q$.
If $|a|\neq 1$, the curves intersect $Q$ only at $(d,d,0)\in\Delta_0$ and $\sigma(d,d,0)\in\Delta_{\infty}$.
The curve \eqref{eqn-001} is always disjoint from $Q$.
(iv) By the birational correspondence between $X$ and $Z^+$, all the curves \eqref{eqn-001} are transformed into  real lines in $Z^+$,
while the curves \eqref{eqn-111} and \eqref{eqn-111'} are transformed into  real curves whose normal bundles are $\mathscr O(2)\oplus \mathscr O$, as long as $|a|=1$.
\end{prop}

The final assertion of (iv) means that if $|a|=1$ the curves \eqref{eqn-111} and \eqref{eqn-111'} are not transformed into  lines.
On the other hand, it will be shown (in the proof of Theorem \ref{thm-main})  that if $|a|\neq 1$ the curves \eqref{eqn-111} and \eqref{eqn-111'} are transformed into  real lines in $Z^+$.

\bigskip
\noindent Proof of Proposition \ref{prop-explicit}.
(i), (ii) and (iii) are elementary and we omit proofs.
For (iv), let $C\subset X$ be a curve defined by \eqref{eqn-001} (for some $d\in\mathbf C\cup\{\infty\}$).
Then we can readily find divisors $D_1$ and $D_2$ in $Z^+$ intersecting transversally along the image of $C$ in $Z^+$, to conclude that $C$ is transformed into a smooth rational curve satisfying $N\simeq\mathscr O(1)\oplus\mathscr O(1)$.
On the other hand, if $C$ is a curve defined by \eqref{eqn-111} or \eqref{eqn-111'} with $|a|=1$, $C$ is contained in $Q$ as in (iii).
It is readily  seen that $N_{C/Q}\simeq\mathscr O(2)$.
If we write $C^+$ for the image of $C$ in $Z^+$, we obtain from the inclusions $C^+\subset Q^+\subset Z^+$ an exact sequence
$$
0\lra \mathscr O(2)\lra N_{C^+/Z^+}\lra N_{Q^+/Z^+}|_{C^+}\lra 0.
$$
Further, there is a natural isomorphism $Q^+\simeq Q\simeq\mathbf P\times\mathbf P$ and we have $N_{Q^+/Z^+}\simeq \mathscr O(2,-2)$, while $C^+$ is a $(1,1)$-curve in $Q^+$.
Hence the restriction $N_{Q^+/Z^+}|_{C^+}$ is trivial and we obtain $N_{C^+/Z^+}\simeq  \mathscr O(2)\oplus\mathscr O$.
\proofend

\bigskip
Since the curve (\ref{eqn-111}) is determined by two numbers  $d\in\mathbf C$ and $a\in\mathbf C^*$, we denote it by  $L_{d,a}$.
Also for $d=\infty$ and $a\in\mathbf C^*$, we denote $L_{\infty,a}$ for the curve \eqref{eqn-111'}, by abuse of notation.
Of course, $d\in\mathbf C\cup\{\infty\}$ specifies the intersection point of the curve with $\Delta_0$.
Further for  $d\in\mathbf C\cup\{\infty\}$ we write $L_{d,0}$ for the curve \eqref{eqn-001} whose degree is $(0,0,1)$.
Then as $a$ goes to $0$ or $\infty$, $L_{d,a}$ degenerates into a reducible curve as follows:

\begin{lemma} \label{lemma-red}
(i)
When the parameter $a$ goes to zero, the curve $L_{d,a}$  degenerates into a reducible, connected curve $L_{d,0}+C_0+C_{\infty}$, where $C_0$ (resp. $C_{\infty}$) is the unique holomorphic curve of degree $(1,0,0)$ (resp.\,$(0,1,0)$) in $X$ going through the point $(d,d,0)$ (resp. \!\!$(-1/\ol{d},-1/\ol{d},\infty)$). 
(ii) When the parameter $a$ goes to  infinity, the curve $L_{d,a}$  degenerates into a reducible, connected curve $L_{-1/\ol{d},0}+C'_0+C'_{\infty}$, where $C'_0$ (resp. $C'_{\infty}$) is the unique holomorphic curve of degree $(0,1,0)$ (resp.\,$(1,0,0)$) in $X$ going through the point $(d,d,0)$  (resp. \!\!$(-1/\ol{d},-1/\ol{d},\infty)$). 
\end{lemma}

\begin{figure}
\includegraphics{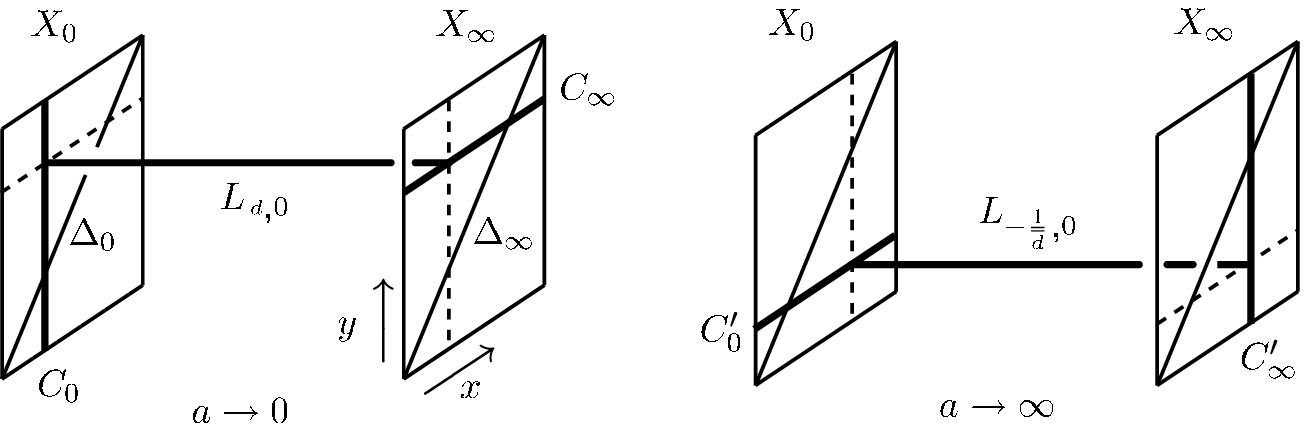}
\caption{The limit curves (bold ones) as $a\to0$ and $a\to\infty$}
\label{fig-limit}
\end{figure}

The proof is also straightforward from the expressions \eqref{eqn-111}--\eqref{eqn-001}, if one keeps in mind that the degree of a curve is  preserved even after taking limits.
An important point in the lemma is that, due to the difference of degrees of the extra components ($C_0$ of degree $(1,0,0)$ and $C'_0$ of degree $(0,1,0)$, and also $C_{\infty}$ of degree $(0,1,0)$ and $C'_{\infty}$ of degree $(1,0,0)$; see Figure \ref{fig-limit}), we have
\begin{equation}\label{diflimit}
\lim_{a\ra 0} L_{d,a}\neq \lim_{a\ra\infty} L_{d',a}
\end{equation}
 for any $d,d'\in\mathbf C\cup\{\infty\}$.

Thus we have obtained a  family of real curves of degree $(1,1,1)$ in $X$.
It consists of the following three kinds of curves:

\vspace{1mm}
\noindent
(A) irreducible curves $L_{d,a}$ defined by \eqref{eqn-111} and \eqref{eqn-111'}, where $d\in\mathbf C\cup\{\infty\}$ and $a\in\mathbf C^*$ with $|a|\neq 1$,

\noindent
(B) irreducible curves $L_{d,a}$ defined by \eqref{eqn-111} and \eqref{eqn-111'}, where $d\in\mathbf C\cup\{\infty\}$ and $|a|=1$

\noindent
(C) reducible curves $L_{d,0}+C_0+C_{\infty}$ as in (i) of Lemma \ref{lemma-red}, and $L_{d',0}+C'_0+C'_{\infty}$ as in (ii) of the lemma, where $d, d'\in\mathbf C\cup\{\infty\}$.

\vspace{1mm}
For  explanation we divide (A) into two subfamilies (A1) and (A2), where (A1) consists of $L_{d,a}$ with $|a|<1$ and (A2) consists of $L_{d,a}$ with $|a|>1$.
Similarly  we define (C1) to be the family consisting of the former reducible curves $L_{d,0}+C_0+C_{\infty}$ in (C), and (C2) to be the family consisting of the latter reducible curves $L_{d',0}+C'_0+C'_{\infty}$ in (C).
Then by Lemma \ref{lemma-red}, the union (A1) $\cup$ (C1) constitutes a connected family and (A2) $\cup$ (C2) constitutes (another) connected  family.
Note thanks to \eqref{diflimit} these two connected families (A1) $\cup$ (C1) and (A2) $\cup$ (C2) contain no common curves.
Let $M^+$ and $M^-$ be the parameter spaces of the former and the latter families of curves respectively.
Let $K$ be the parameter space of the curves (B).
Since  curves in (A) and (B) are clearly deformation equivalent, $M^+$ and $M^-$ are connected by $K$ and we obtain a compact, connected 4-manifold 
\begin{equation}\label{dcp3}
M:=M^+\cup K\cup M^-.
\end{equation}
$M$ parameterizes all curves in (A), (B) and (C).
By the projection assigning $d\in \mathbf C\cup\{\infty\}=\mathbf P$, $M^+$ and $M^-$ have structures of  disk bundles over $\mathbf P$ and $K$ has a structure of a circle bundle over $\mathbf P$.
Consequently, $M$ has a structure of an $S^2=\mathbf P$ bundle over $\mathbf P$.

\begin{rmk}
{\em
The 4-manifold $M$ has a natural structure of a differential manifold which has $(d,a)$, $(1/d,a)$, $(d,1/a)$ and $(1/d,1/a)$ as complex local coordinates.
From algebraic geometric point of view, $M$ is a connected component of a real locus of a Hilbert scheme of curves in  $X=\mathbf P\times\mathbf P\times\mathbf P$ whose degrees are $(1,1,1)$ and which intersect both of $\Delta_0$ and $\Delta_{\infty}$.
(This Hilbert scheme is readily shown to be 4-dimensional (over $\mathbf C$).)
}
\end{rmk}

The involution of $X$ interchanging the first and the second factors naturally induces an involution  of $M$. It is explicitly given by 
 \begin{equation}\label{eqn-glue}
(d,a)\mapsto \left(d,\,\frac{1}{\ol{a}}\right).
\end{equation}
and it follows that the involution of $M$ fixes every points of $K$ (which is obvious from the construction) and induces an isomorphism between $M^+$ and $M^-$.

Then the following lemma is important in proving our main result:

\begin{lemma}\label{lemma-map}
Let $M$, $K$ and $M\backslash K=M^+\cup M^-$ be as above, and fix $t\in C$ with $t\neq 0,\infty$. 
Let  $u:M\ra X_t=\mathbf{P}\times\mathbf{P}$ be a (non-holomorphic) map defined by
$$u(d,a)=\left(\frac{d-at}{1+a\ol{d}t},\,\,\frac{\ol{a}d-t}{\ol{a}+\ol{d}t}\right).$$
(Namely $u$ assigns the intersection point $L_{d,a}\cap X_t$ for each $(d,a)\in M$.)
Then 
(i) the map $u$ is surjective,
(ii) $u(K)=\Delta_t$ (= the diagonal of $X_t$) and $u^{-1}(\Delta_t)=K$,
(iii) $u$ maps $M^+$ and $M^-$ diffeomorphically onto $X_t\backslash\,\Delta_t$.
\end{lemma}

\noindent Proof.
First we show  $u$ is surjective. 
By eliminating $a$ from (\ref{eqn-111}), we readily see that, if we fix $d\in\mathbf C$ and $t\in C$ ($t\neq 0,\infty$),
then  the set of intersection points $\{L_{d,a}\cap X_t\set a\in\mathbf C^*\}$ is contained in an antiholomorphic curve $\mathscr C_{d,t}$ in $X_t$ defined by 
\begin{equation}\label{eqn-traject}
\mathscr C_{d,t}:\,\,y=\frac{d(1+R)\ol{x}+(R-|d|^2)}{(|d|^2R-1)\ol{x}+\ol{d}(1+R)},
\end{equation}
where we set $R=|t|^2$ for simplicity, and $(x,y)$ is an affine coordinate on $X_t=\mathbf P\times\mathbf P$ as before.
Note that (\ref{eqn-traject}) is the graph of (an anti-holomorphic) fractional transformation.
It is easy to see that $\mathscr C_{d,t}$ is not a graph of a constant function. 
On the other hand,  viewing the $x$-factor of $u$ as a fractional transformation of the variable $a$, the map $\mathbf P\ni a\mapsto x(u(d,a))\in\mathbf P$ is readily verified to be  surjective.
It follows that the map $\mathbf P\ni a\mapsto u(d,a)\in\mathscr C_{d,t}$ is surjective 
for any fixed $d\in\mathbf P$ and $t\neq 0,\infty$.
Therefore, in order to prove the surjectivity of the map $u:M\ra X_t$, it suffices to show that the map $d\mapsto y$ defined by (\ref{eqn-traject}) is surjective, while $x$ and $t$ are fixed.
To see this, we introduce a new variable $c$ by setting
$$
c=\frac{d-x}{\ol{x}d+1},
$$
which is a non-degenerate fractional transformation of $d$.
Then the equation of $\mathscr C_{d,t}$ becomes
\begin{equation}\label{eqn-traject'}
y=-\frac{(1-\ol{c}x)R-\ol{c}\,(c+x)}{(\ol{c}+\ol{x})R+\ol{c}\,(1-c\ol{x})}.
\end{equation}
It suffices to show that the map $c\mapsto y$ defined by (\ref{eqn-traject'}) is surjective.
It is readily seen that (\ref{eqn-traject'}) is the composition of two maps
$$
c\mapsto \frac{R}{\ol{c}}-c\,(=:b) \hspace{3mm}\mbox{and}\hspace{3mm}
 b\mapsto -\frac{b-x(1+R)}{\ol{x}{b}+(1+R)}.
$$
The former map contracts the circle $|c|=\sqrt{R}$ to a point $\{0\}$ and maps the remaining two discs diffeomorphically onto $\mathbf P\bsl \{0\}$. 
The latter fractional transformation is readily verified to be biholomorphic on $\mathbf P$, which maps the point $\{b=0\}$ to the point $\{x\}$.
Thus, the map $c\mapsto y$ is surjective, contracting the circle $|c|=|t|$ to a point $\{x\}$.
Equivalently,  $d\mapsto y$ is surjective, contracting the circle $|x-d|/|\ol{d}x+1|=|t|$ to a point $\{x\}$.
In particular we obtain (i) of the lemma.

(ii) is immediate from Proposition \ref{prop-explicit} (iii).

Finally we prove (iii).
Since the Jacobian of $u$  is calculated to be
$$
\frac{(1+|d|^2)(1+|t|^2)|t|^2}{|1+a\ol{d}t|^4|a+d\ol{t}|^4}\left(|a|^4-1\right),
$$ 
$u$ can fail to be locally diffeomorphic only on the set $K=\{|a|=1\}$.
Hence $u:M\bsl K\ra\mathbf P\times\mathbf P\bsl\Delta$ is an unramified double covering.
Moreover since $\mathbf P\times\mathbf P\bsl\Delta$ is simply connected, $u:M\bsl K=M^+\cup M^-\ra\mathbf P\times\mathbf P\bsl\Delta$ must be diffeomorphic on $M^+$ and $M^-$ respectively.
This proves (iii) of the lemma.
\proofend

\bigskip

So far we have investigated a family of curves in $X$ that can be transformed into real lines in $Z^+$. 
Their parameter spaces are 4-manifolds $M^+$ and $M^-$.
Now we consider images of the curves under the birational transformation from $X$ to $Z^+$.
The results are as follows:

\begin{itemize}
\item
All curves  belonging to (A) and (B) are transformed into real irreducible curves in $Z^+$.
(Recall that curves in (B) are not transformed into   lines by Proposition \ref{prop-explicit} (iv).)
\item
Reducible curves belonging to (C1) are transformed into real {\em irreducible}\, curves  in $Z^+$. (The two extra components $C_0$ and $C_{\infty}$ are contracted to points by $\nu^+:Y\ra Z^+$. See Figure \ref{fig-limit}).
\item
For reducible curves $L_{d,0}+C'_0+C'_{\infty}$ belonging to (C2), $L_{d,0}$ is transformed into real line as in Proposition \ref{prop-explicit} (iv), and $C'_0$ and $C'_{\infty}$ are transformed into $B^+_0$ and $B^+_{\infty}$ respectively.
(Note that $B_0^++B_{\infty}^+$ is homologous to zero as mentioned in the first paragraph of this section.)
\end{itemize}

We remark  that if we consider $\nu^-:Y\ra Z^-$ instead of $\nu^+:Y\ra Z^+$, the role of (C1) and (C2) in the final item are interchanged.

Let $\mathscr L^+$ (resp. $\mathscr L^-$) be the family of curves in $Z^+$ whose members are the transformations of curves in $X$ parameterized by $M^+$ (resp. $M^-$).
Then the above 3 items mean that all members of $\mathscr L^+$ are irreducible, while $\mathscr L^-$ contains a 2-dimensional family of reducible members.
Moreover, members of $\mathscr L^+$ and $\mathscr L^-$ are connected by a family of real curves which are not lines (parameterized by $K$).
In particular, they are homologous in $Z^+$.
Then the following is the main result:

\begin{thm}\label{thm-main}
Let $Z^+$ and \,$Q^+\subset Z^+$  be as before, and $\mathscr L^+$ the family of real irreducible curves in $Z^+$ as above.
Then the complement $Z^+\backslash Q^+$ has a structure of the twistor space of a hyperK\"ahler metric on the cotangent bundle of\, $\mathbf P$, which has $\mathscr L^+$ as the set of twistor lines. 
\end{thm}

By a result of Fujiki \cite{F} the hyperK\"ahler metric  we obtain is precisely the Eguchi-Hanson metric \cite{EH78}. Note that another family $\mathscr L^-$ in $Z^+$ cannot be the set of twistor lines since it contains reducible members.
 
 \vspace{3mm}
\noindent Proof of Theorem \ref{thm-main}.
According to \cite[Theorem 3.3]{HKLR},  it suffices to show 
(a) $Z^+\backslash Q^+$ is foliated by  members of $\mathscr L^+$;  namely $Z^+\backslash Q^+$ is covered by curves of $ \mathscr L^+$, and that different curves of $\mathscr L^+$ are always disjoint,
(b) there is a non-vanishing holomorphic section of $\wedge^2 T^*_{Z/C}\otimes f_+^*\mathscr O(2)$, 
(c) the parameter space $M^+$ of $\mathscr L^+$ is diffeomorphic to the cotangent bundle of $\mathbf P$.
(Theorem 3.3 of  \cite{HKLR} requires that the normal bundles of the curves in $\mathscr L^+$ are isomorphic to $\mathscr O(1)\oplus \mathscr O(1)$.
But this follows from (a). In fact it is easy to show that the normal bundles in question are isomorphic to $\mathscr O(1)\oplus\mathscr O(1)$ or $\mathscr O(2)\oplus\mathscr O$.
If there is a curve such that the latter holds, there has to exist a real 3-dimensional subfamily of $\mathscr L^-$ such that every two members of the subfamily intersect at two points. This contradicts (a).)

We first verify (a).
It is clear that $Z^+\backslash (Z_0\cup Z_{\infty})$ is naturally biholomorphic to $X\backslash( X_0\cup X_{\infty})$ ($\simeq \mathbf P\times\mathbf P\times\mathbf C^*$), where $Z_t=f_+^{-1}(t)$ and $X_t=f^{-1}(t)$ for $t\in C$.
So any point $z$ on $Z^+\backslash (Z_0\cup Z_{\infty})$ can be specified by $(x,y,t)$ for some $t\neq 0,\infty$ and $(x,y)\in \mathbf P\times\mathbf P$.
We have $(Z^+\bsl Q^+)\cap Z_t\simeq X_t\bsl\Delta_t$.
Therefore, Lemma \ref{lemma-map} directly implies that $z=(x,y,t)\in Z^+\bsl Q^+$  is passed by a unique member of $\mathscr L^+$. 
Next take any $z\in Z_0\backslash Q^+$.
By construction  $Z_0$ is identified with the projectified normal bundle $\mathbf P(N_{\Delta_0/X})$ and we have $N_{\Delta_0/X}\simeq N_{\Delta_0/X_0}\oplus N_{\Delta_0/Q}$, where $Q=\Delta\times C\subset X$ as before.
Thus the intersection point of a member of $\mathscr L^+$ with $Z_0$ is represented by a tangent vector of the corresponding curve in $X$ at its intersection point with $X_0$.
If $d\neq \infty$, by (\ref{eqn-111}), a tangent vector $v$ at the point is readily seen to be
\begin{equation}
v=\left(-a\left(1+|d|^2\right),\,-\frac{1}{\ol{a}}\left(1+|d|^2\right),\,1\right).
\end{equation}
On the other hand, because $\Delta_0$ is defined by $x-y=t=0$ in $X$, we can take $(x-y,t)$ as a linear coordinate on the vector space $T_{\xi}X/T_{\xi}\Delta_0\,(\simeq\mathbf C^2)$, where $\xi=(d,d,0)$ denotes the intersection point.
Thus the point $Z_0\cap L_{d,a}$ is represented by a vector 
\begin{equation}\label{eqn-2:1'}
\left((1+|d|^2)\left(a-\frac{1}{\ol{a}}\right), 1\right).
\end{equation}
It is easy to see that a map from $\mathbf C^*$ to $\mathbf C$ defined by 
$a\mapsto a-(1/\ol{a})$
has the following properties: 
it contracts $U(1)=\{ |a|=1\}$ to the origin, and maps each connected component of $\mathbf C^*\bsl U(1)$ diffeomorphically onto $\mathbf C^*$.
By (\ref{eqn-2:1'}) it follows that if we fix $d\in\mathbf C$ and allow the parameter  $a$ to move in the domain $|a|<1$,  then the intersection point with $Z_0$ takes each point of  the fiber of $Z_0\ra\Delta_0$ over $\xi$ precisely once, except the two points corresponding to the two vectors represented by $(1,0)$ and $(0,1)$ in the above coordinate.
The point corresponding to $(0,1)$ (which is represented by a trivial section of $f:X\ra C$) is then passed by the transformation of a curve of (C1).
On the other hand, the point corresponding to $(1,0)$ (which is represented by a vector tangent to $X_0=f^{-1}(0)$) is contained in $Q^+$.
So the point does not need to be passed.
It remains to show that the fiber of $Z_0\ra \Delta_0$ (outside $Q^+$) over $(\infty,\infty,0)\in\Delta_0$ is passed by some $L_{{\infty},a}$ for a unique $a$ with $|a|<1$.
But this can be also verified by easier calculation using (\ref{eqn-111'}) instead of (\ref{eqn-111}).
Furthermore, all points of $Z_{\infty}\backslash Q^+$ are also passed by a unique member of $\mathscr L^+$, since $Z_{\infty}=\sigma(Z_0)$ and all members of $\mathscr L^+$ are real.
Thus we have shown that $Z^+\bsl Q^+$ is foliated by members of  $\mathscr L^+$, proving the assertion (a).

Next  we show the assertion (b).
Since $f_+:Z^+\ra C$ is a holomorphic submersion, there is  an exact sequence
\begin{equation}
0\lra T_{Z^+/C}\lra T_{Z^+}\lra f^*T_C\lra 0.
\end{equation}
It follows that $\det T_{Z^+}\simeq \wedge^2 T_{Z^+/C}\otimes f^*\mathscr O(2)$, and equivalently, $K_{Z^+}\simeq\wedge^2 T^*_{Z^+/C}\otimes f^*\mathscr O(-2)$.
On the other hand, by (\ref{eqn-can}), we have $K_{Z^+}\simeq \mathscr O(-2Q^+)\otimes f_+^*\mathscr O(-4)$.
From these we obtain $\wedge^2T^*_{Z^+/C}\simeq  \mathscr O(-2Q^+)\otimes f^*\mathscr O(-2)$.
From this we conclude $\wedge^2T^*_{Z^+/C}\otimes f^*\mathscr O(2)$ is trivial over $Q^+$.
Thus we have shown (b).

(c) is obvious since we have just seen that any point of $Z_0\backslash Q^+$ is passed by a unique member of $\mathscr L^+$ and since $Z_0\backslash Q^+$ is isomorphic to the total space of $\mathscr O(2)$.
This finishes a proof of Theorem \ref{thm-main}.
\proofend

\bigskip

As in $Z^+$, there are two families of real curves in $Z^-$ whose parameter spaces are $M^+$ and $M^-$.
If we again write them by $\mathscr L^+$ and $\mathscr L^-$ respectively, the complement $Z^-\backslash Q^-$ has a structure of the twistor space of a hyperK\"ahler metric having $\mathscr L^-$ as the set of twistor lines.
Of course this is isomorphic to the hyperK\"ahler metric obtained in Theorem \ref{thm-main}.

Finally we give some remarks on the automorphism group of Nagata threefold and its induced action on the 4-manifold $M=M^+\cup K\cup M^-$.
The identity component of the holomorphic automorphism group Aut$_0(X)$ of $X=\mathbf P\times\mathbf P\times\mathbf P$ is of course the product  of three copies of the Lie group $PSL(2,\mathbf C)$.
An automorphism $g=(g_1,g_2,g_3)\in$\, Aut$_0(X)$ ($g_i\in PSL(2,\mathbf C)$) induces that on the Nagata threefold iff $g$ preseves $\Delta_0$ and $\Delta_{\infty}$ (which are the center of the blowing-up $\mu:Y\to X$).
This condition is equivalent to the conditions $g_1=g_2$ and $g_3\in\mathbf C^*$, where $\mathbf C^*$ is the subgroup consisting of automorphisms which fix $0$ and $\infty\in\mathbf P$.
Conversely, if an automorphism $h$ of the Nagata threefold belongs to an identity component of the holomorphic automorphism group, 
$h$ preserves the center of the blowing-down $\nu^{\pm}:Y\to Z^{\pm}$.
Hence $h$ always induces an automorphism on $X$, which necessarily belongs to the identity component.
These mean that  the identity component of the holomorphic automorphism group of the Nagata threefold is the subgroup $PSL(2,\mathbf C)\times\mathbf C^*$ of Aut$_0(X)$ given by the injection $(g_1,g_3)\mapsto (g_1,g_1,g_3)$.
Further, an automorphism $g=(g_1,g_1,g_3)$ ($g_1\in PSL(2,\mathbf C)$ and $g_3\in\mathbf C^*$) commutes with the real structure $\sigma$ (defined by \eqref{eqn-rs}) 
iff $g_1\in PSU(2)$ and $g_3\in U(1)$.
This means that the identity component of  the group of holomorphic automorphisms of the Nagata threefold commuting with the real structure is exactly $PSU(2)\times U(1)$.

Since the action of $PSU(2)\times U(1)$ preserves $\Delta_0\cup\Delta_{\infty}$ and commutes with the real structure, it naturally operates on the 4-manifold $M$.
If $g_1\in PSU(2)$ is represented by the matrix
$$
\begin{pmatrix}
\alpha & \beta\\
-\ol{\beta} &\ol{\alpha}
\end{pmatrix}
$$
satisfying $|\alpha|^2+|\beta|^2=1$,
then $(g_1,g_3)\in PSU(2)\times U(1)$ maps the curve $L_{d,a}$ (defined by the equation \eqref{eqn-111}) to the curve $L_{d',a'}$ where $d'$ and $a'$ satisfy
\begin{equation}\label{tra}
d'=\frac{\alpha d+\beta}{-\ol{\beta}d+\ol{\alpha}},
\,\,\,a'=\frac{\alpha-\beta\ol{d}}{\ol{\alpha}-\ol{\beta} d}\cdot g_3a
\end{equation}
This shows that the $PSU(2)\times U(1)$-action on $M$ preserves the decomposition $M=M^+\cup K\cup M^-$ obtained in \eqref{dcp3}.
Also it follows from \eqref{tra} that the $PSU(2)\times U(1)$-action on the $U(1)$-bundle $K$ (over $\mathbf P$) is transitive and its isotropy subgroup at the point $(d,a)=(0,1)$ is $\{{\rm{id}}\}\times U(1)$.
Thus $K$ is diffeomorphic to $PSU(2)=SO(3)$.
Further, the $U(1)$-bundle map $K\to\mathbf P$ (assigning $d$) is exactly obtained from the Hopf fibration $S^3\to S^2$ by dividing each fiber by multiplying $-1\in U(1)$.
This means that $K$ is identified with the unit circle bundle of the line bundle $\mathscr O(-2)\to\mathbf P$.
Therefore, the disk bundles $M^+\to\mathbf P$ and $M^-\to \mathbf P$ (also assigning $d$) are identified with a disk bundle of $\mathscr O(-2)$.

\small
\vspace{13mm}
\hspace{7.5cm}
$\begin{array}{l}
\mbox{Department of Mathematics}\\
\mbox{Graduate School of Science and Engineering}\\
\mbox{Tokyo Institute of Technology}\\
\mbox{2-12-1, O-okayama, Meguro, 152-8551, JAPAN}\\
\mbox{{\tt {honda@math.titech.ac.jp}}}
\end{array}$


\begin{thebibliography}{99}


\bibitem{ahs}M. Atiyah, N. Hitchin, I. Singer,
 {\em  Self-duality in four-dimensional
Riemannian geometry},
Proc. Roy. Soc. London, Ser. A {\bf 362} (1978),
425--461.

\bibitem{B98}
O. Biquard, {\em Twisteurs des orbites coadjointes et m\'etriques hyper-pseudo K\"ahl\'eriennes},
Bull. Soc. Math. France {\bf 126} (1998), 79--105.

\bibitem{B86}
D. Burns,
{\em Some examples of the twistor constructions},
in ``Contributions to several complex variables: in honor of Wilhelm Stoll''
(A. Howard and P. -M. Wong ed.)
Vieveg, Braunschweig (1986), 51--67.

\bibitem{EH78} T. Eguchi, A.J. Hanson,
{\em Asymptotically flat solutions to Euclidean gravity},
Phys. Lett. {\bf 74B} (1978), 249--251.

\bibitem{F} A. Fujiki,
{\em  Nagata threefold and twistor space},
Quaternionic structures in mathematics and physics (Trieste, 1994), 139--146 (electronic), Int. Sch. Adv. Stud. (SISSA), Trieste, 1998. 

\bibitem{Hi79}N. Hitchin, {\em Polygons and Gravitons},
Math. Proc. Cambridge Philos. Soc. {\bf 85} (1979),  465--476.

\bibitem{HKLR}
N. J. Hitchin,  A. Karlhede, U. Lindstr\"om, M. Ro\v cek, 
{\em Hyper-K\"ahler metrics and supersymmetry},
Comm. Math. Phys. {\bf 108} (1987), 535--589.


\bibitem{Na58}M. Nagata, 
{\em Existence theorems for non-projective
complete algebraic varieties},
Ill. J. Math {\bf 2} (1958),  490--498.



\end{thebibliography}
\end{document}